\documentclass[12pt]{amsart}

\usepackage{amsfonts,amsmath,amsthm}
\usepackage{latexsym}

\newtheorem{theorem}{Theorem}[section]

\newtheorem{lemma}[theorem]{Lemma}
\newtheorem{remark}[theorem]{Remark}

\theoremstyle{definition}

\newcommand{\F}{\ensuremath{\mathbb{F}}}

\newcommand{\ZZ}{\ensuremath{\mathbb{Z}}}

\def \ffi {\varphi}

\def \< {\langle}
\def \> {\rangle}

\begin{document}

\title{Squares and difference sets in finite fields}

\author[C. Bachoc]{C. Bachoc}
\address{C. B.: Univ Bordeaux,
Institut de Math\'ematiques de Bordeaux, 351,
cours de la Lib\'eration 33405, Talence cedex, France,
Tel: (+33) 05 40 00 21 61, Fax: (+33) 05 40 00 21 23}
\email{bachoc@math.u-bordeaux1.fr}

\author[M. Matolcsi]{M. Matolcsi}
\address{M. M.: Alf\'ed R\'enyi Institute of Mathematics,
Hungarian Academy of Sciences POB 127 H-1364 Budapest, Hungary
Tel: (+361) 483-8307, Fax: (+361) 483-8333}
\email{matomate@renyi.hu}

\author[I. Z. Ruzsa]{I. Z. Ruzsa}
\address{I. Z. R.: Alfr\'ed R\'enyi Institute of Mathematics,
Hungarian Academy of Sciences POB 127 H-1364 Budapest, Hungary
Tel: (+361) 483-8328, Fax: (+361) 483-8333}
\email{ruzsa@renyi.hu}

\thanks{M.M. and I.Z.R. were supported by the ERC-AdG 228005,
and OTKA Grants No. K81658, and M.M. also by the Bolyai
Scholarship.}

\begin{abstract}
For infinitely many primes $p=4k+1$ we give a slightly improved upper bound for the maximal cardinality of a set $B\subset \ZZ_p$ such that the difference set $B-B$ contains only quadratic residues. Namely, instead of the "trivial" bound $|B|\leq \sqrt{p}$ we prove $|B|\leq \sqrt{p}-1$, under suitable conditions on $p$. The new bound is valid for approximately three quarters of the primes $p=4k+1$.

\bigskip

\end{abstract}

\maketitle

\bigskip


\bigskip

\section{introduction}

Let $q$ be a prime-power, say $q=p^k$. We will be interested in estimating the maximal cardinality $s(q)$
of a set $B\subset \F_q$ such that the difference set $B-B$ contains only squares. While our main interest is in the case $k=1$,
we find it instructive to compare the situation for different values of $k$.

This problem makes sense only if $-1$ is a square; to ensure this we assume $q\equiv 1 \pmod 4$.
The universal upper bound $s(q)\leq \sqrt{q}$ can be proved by a pigeonhole argument or by simple Fourier anlysis, and it
 has been re-discovered several times (see \cite[Theorem 3.9]{delsarte}, \cite[Problem 13.13]{lovasz}, \cite[Proposition 4.7]{cameron}, \cite[Chapter XIII, Theorem 14]{bollobas}, \cite[Theorem 31.3]{lint}, \cite[Proposition 4.5]{sudakov}, \cite[Section 2.8]{lev}
for various proofs).
 For even $k$ we have equality, since $\F_{p^k}$ can be constructed as a quadratic extension of $\F_{p^{k/2}}$, and then
every element of the embedded field $\F_{p^{k/2}}$ will be a square. It is known that every case of equality can be obtained by a linear
transformation from this one, \cite{blokhuis}.


Such problems and results are often formulated in terms of
the Paley graph $P_q$, which is the graph with vertex set $\F_q$ and an edge between $x$ and $y$ if and only if $x-y =a^2$ for some non-zero $a \in \F_q$.

Paley graphs are self-complementary, vertex and
edge transitive, and $(q, (q-1)/2, (q-5)/4, (q-1)/4)$-strongly regular (see \cite{bollobas} for these and other basic properties of $P_q$). Paley graphs have received considerable attention over the past decades because they exhibit many properties of random graphs $G(q, 1/2)$ where each edge is present with probability $1/2$. Indeed, $P_q$ form a family of {\it quasi-random} graphs, as shown in \cite{quasi}.

With this terminology $s(q)$ is the {\it clique number} of $P_q$. The general lower bound $s(q)\geq (\frac{1}{2}+o(1)) \log_2 q$ is established in \cite{cohen}, while it is proved in \cite{graham}
that $s(p)\geq c \log p \log \log \log p$ for infinitely many primes $p$. The ``trivial'' upper bound $s(p)\leq \sqrt{p}$
is notoriously difficult to improve, and it is  mentioned explicitly in the selected list of problems \cite{lev}. The only improvement we are aware of concerns the special case $p=n^2+1$ for which it is proved in \cite{penman} that $s(p)\leq n-1$ (the same result was proved independently by T. Sanders -- unpublished, personal communication). It is more likely, heuristically, that the lower bound is closer to the truth than the upper bound. Numerical data \cite{shearer, exoo} up to $p<10000$ suggest (very tentatively) that the correct order of magnitude for the clique number of $P_p$ is $c \log^2 p$ (see the discussion and the plot of the function $s(p)$ at \cite{web}).

\medskip

In this note we prove the slightly improved upper bound $s(p)\leq \sqrt{p}-1$ for the {\it majority} of the primes $p=4k+1$ (we will often suppress the dependence on $p$, and just write $s$ instead of $s(p)$).



\medskip
We will denote the set of nonzero quadratic residues by $Q$, and that of nonzero non-residues by $NQ$. Note that $0\notin Q$ and $0\notin NQ$.

\section{The improved upper bound}

\begin{theorem}\label{newbound}
Let $q$ be a prime-power, $q=p^k$, and assume that $k$ is odd and $q\equiv 1 \pmod 4$. Let  $s=s(q)$ be the maximal cardinality
of a set $B\subset \F_q$ such that the difference set $B-B$ contains only squares. \\
(i) If $[\sqrt{q}]$ is even then $s^2+s-1\leq q$,\\
(ii) if $[\sqrt{q}]$ is odd then $s^2+2s-2\leq q$.
\end{theorem}

\begin{proof}
The claims hold if $s<[\sqrt{q}]$. Hence we may assume that $s\geq [\sqrt{q}]$.

\begin{lemma} Let $D\subset\F_q$ be a set such that
   \[ D\subset NQ, \ D-D \subset Q\cup \{0\}. \]
With $r=|D|$ we have
\begin{equation}\label{sbound}
s(q)\leq 1+\frac{q-1}{2r}.
\end{equation}
\end{lemma}

\begin{proof}
  Let $B$ be  a maximal set such that $B-B\subset Q\cup \{0\}$, $|B|=s(q)=s$.
  Consider the equation
   \[ b_1-b_2 = zd , \ b_1, b_2\in B, \ d\in D, \ z\in NQ . \]
   This equation has exactly $s(s-1)r $ solutions; indeed, every pair of distinct $ b_1, b_2\in B$ and a $d\in D$ determines $z$ uniquely.
On the other hand, given $b_1$ and $z$, there can be at most one pair $b_2$ and $d$ to form a solution. Indeed, if there were
another pair $b_2', d'$, then by substracting the equations
 \[  b_1-b_2 = zd , \  b_1-b_2' = zd'  \ \]
 we get $(b_2'-b_2)=z(d-d')$, a contradiction, as the left hand side is a square and the right hand side is not.
This gives $s(s-1)r \leq s(q-1)/2$ as wanted.
\end{proof}

We try to construct such a set $D$ in the form $D=(B-t) \cap NQ$ with a suitable $t$. The required property then follows from
$D-D\subset B-B$.







Let $\chi$ denote the quadratic multiplicative character, i.e. $\chi(t)=\pm 1$ according to whether $t\in Q$ or $t\in NQ$ (and $\chi(0)=0$). Let
\begin{equation}\label{phi}
\ffi(t)=\sum_{b\in B} \chi(b-t).
\end{equation}
Clearly
 \[ \varphi(t) = |(B-t)\cap Q| - |(B-t)\cap NQ|,  \]
and hence for $t \notin B$ we have
  \[  |(B-t)\cap NQ| = \frac{s-\varphi(t)}{2}. \]
To find a large set in this form we need to find a negative value of $\varphi$.



We list some properties of this function. For $t\in B$ we have $\varphi(t)=s-1$, and otherwise
 \[ \varphi(t) \leq s-2, \ \varphi(t)  \equiv s \pmod 2\]
(the inequality expresses the maximality of $B$). Furthermore,
  \[ \sum_{t} \varphi(t)= 0, \]
  and, since translations of the quadratic character have the quasi-orthogonality property
   \[ \sum_t \chi(t+a) \chi(t+b) = -1 \]
   for $ a \neq b$, we conclude
 \[ \sum_{t} \varphi(t)^2 = s(q-1) - s(s-1) = s(q-s) . \]
 By substracting the contribution of $t\in B$ we obtain
\[ \sum_{t\notin B } \varphi(t)= -s(s-1), \]

 \[ \sum_{t \notin B } \varphi(t)^2 = s(q-s) - s(s-1)^2 = s(q-s^2+s-1) . \]

These formulas assume an even nicer form by introducing the function $\varphi_1(t)=\varphi(t)+1$:
\begin{equation}\label{1mom}
 \sum_{t\notin B } \varphi_1(t)= q-s^2,  \end{equation}
\begin{equation}\label{2mom}
  \sum_{t \notin B } \varphi_1(t)^2 = (s+1)(q-s^2).  \end{equation}
As a byproduct, the second equation shows the familiar estimate $s\leq \sqrt{q}$, so we have $s= [\sqrt{q}] <\sqrt{q}$ (recall that we assume that $s\ge [\sqrt{q}]$, the theorem being trivial otherwise).

Now we consider separately the cases of odd and even $s$. If $s$ is even, then, since $ \sum_{t\notin B } \varphi(t)<0$ and each summand is
even, we can find a $t$ with $\varphi(t)\leq -2$. This gives us an $r$ with $r \geq (s+2)/2$, and on substituting this into
\eqref{sbound} we obtain the first case of the theorem.

If $s$ is odd, we claim that there is a $t$ with $\varphi(t)\leq -3$. Otherwise we have $\varphi(t)\geq -1$, that is, $\varphi_1(t)\geq 0$ for all
$t\notin B$. We also know $\varphi(t)\leq s-2$, $\varphi_1(t) \leq s-1$ for $t\notin B$. Consequently
 \[  \sum_{t \notin B } \varphi_1(t)^2 \leq (s-1) \sum_{t\notin B } \varphi_1(t) = (s-1) (q-s^2), \]
a contradiction to \eqref{2mom}. (Observe that to reach a contradiction we need that $q-s^2$ is strictly positive.
In case of an even $k$ it can happen that $q=s^2$ and the function $\varphi_1$ vanishes outside $B$.)

This $t$ provides us with a set $D$ with  $r \geq (s+3)/2$, and on substituting this into
\eqref{sbound} we obtain the second case of the theorem.
\end{proof}

\begin{remark}\rm
An alternative proof for the case $q=p$ and $s$ being odd is as follows. Assume by contradiction that $\varphi_1$ is even-valued and nonnegative. Then by \eqref{1mom} it must be 0 for at least
   \[ q- |B| - \frac{q-s^2}2 = \frac{q+s^2-2s}2 \]
values of $t$.
 Let $\tilde{\chi}, \tilde{\varphi}, \tilde{\varphi}_1$ denote the images of $\chi, \varphi, \varphi_1$ in $\F_q$ (i.e. the functions are evaluated$\mod p$). By the previous observation $\tilde\varphi_1$ has at least $(q+s^2-2s)/2 $ zeroes.
On the other hand, we have $\tilde\chi(x)=x^{\frac{q-1}{2}}$, and hence $\tilde\varphi_1$ is a polynomial of degree $(q-1)/2$; its leading coefficient is $s=[\sqrt{q}]\neq 0$ mod $p$ (This last fact may fail if $q=p^k$, even if $k$ is odd. Therefore this proof is restricted in its generality. Nevertheless we include it here, because we believe that it has the potential to lead to stronger results if $q=p$.)
Consequently $\tilde\varphi_1$ can have at most  $(q-1)/2$ zeros, a contradiction. In the case of even $k$ we can have $s=\sqrt{q}\equiv0 \pmod p$
and so the polynomial $\tilde\varphi_1$ can vanish, as it indeed does when $B$ is a subfield.
\end{remark}

\begin{remark}\rm
It is clear from \eqref{sbound} that any improved lower bound on $r$ will lead to an improved upper bound on $s$. If one thinks of elements of $\ZZ_p$ as being quadratic residues randomly with probability $1/2$, then we expect that $r\geq \frac{s}{2}+ c\sqrt{s}$. This would lead to an estimate $s\leq \sqrt{p}-cp^{1/4}$. This seems to be the limit of this method. In order to get an improved lower bound on $r$ one can try to prove non-trivial upper bounds on the third moment $\sum_{t\in \ZZ_p} \varphi^3(t)$. To do this, we would need that the distribution of numbers $\frac{b_1-b_2}{b_1-b_3}$ is approximately uniform on $Q$ as $b_1, b_2, b_3$ ranges over $B$. This is plausible because if $s\approx \sqrt{p}$ then the distribution of $B-B$ must be close to uniform on $NQ$. However, we could not prove anything rigorous in this direction.
\end{remark}

\begin{remark}\rm
Theorem \ref{newbound} gives the bound $s\leq [\sqrt{p}]-1$ for about three quarters of the primes $p=4k+1$. Indeed, part {\it (ii)} gives this bound for almost all $p$ such that $n=[\sqrt{p}]$ is odd, with the only exception when $p=(n+1)^2-3$. Part {\it (i)} gives the improved bound $s\leq n-1$ if $n^2+n-1> p$. This happens for about half of the primes $p$ such that $n$ is even.
\end{remark}

\begin{center}Acknowledgment
\end{center}
The authors are grateful to P\'eter Csikv\'ari for insightful comments regarding the prime-power case.

\end{document}